\documentclass[10pt]{amsart}
\usepackage{latexsym}
\usepackage{amssymb}
\usepackage{amsmath}
\usepackage{mathrsfs}
\usepackage{graphicx}
\usepackage{bbm}
\usepackage{bbold}
\input epsf
\input xy
\xyoption{all}
\input diagxy

\def\Aff{\operatorname{Aff}}

\def\co{\operatorname{co}}

\def\spt{\operatorname{spt}}

\begin{document}
\title{Alexandrov's Approach to the Minkowski Problem}
\author{S.~S. Kutateladze}

\address[]{
Sobolev Institute of Mathematics\newline
\indent 4 Koptyug Avenue\newline
\indent Novosibirsk, 630090
\indent RUSSIA}
\email{
sskut@member.ams.org
}
\begin{abstract}
This  article is dedicated to the centenary of the birth
of Aleksandr D. Alexandrov (1912--1999). His
functional-analytical approach to the solving of the
Minkowski problem is examined and applied to the  extremal  problems of
isoperimetric  type with conflicting goals.
\end{abstract}
\date{August 30, 2012}

\maketitle


The {\it Mathematics Subject Classification}, produced jointly by the editorial staffs
of {\it Mathematical Reviews} and {\it Zentralblatt f\"ur Mathematik} in 2010, has Section 53C45
``Global surface theory (convex surfaces \`a la A.~D. Aleksandrov).''
This article surveys some mathematics of the sort.

Good mathematics starts as a first love. If great, it turns into adult sex
and happy marriage. If ordinary, it ends in dumping, cheating or divorce.
If awesome, it becomes eternal.
Alexandrov's mathematics is great (see \cite{SW-I}-\cite{SW-II}). To demonstrate, inspect his solution of the Minkowski problem.

Alexandrov's mathematics is alive, expanding and flourishing for decades.
Dido's problem in the today's setting is one of the examples.

\section*{The Space of Convex Bodies}

A~{\it convex figure\/} is a~compact convex set. A~{\it convex body\/}
is a~solid convex figure.
The {\it Minkowski  duality\/} identifies
a~convex figure $S$ in
$\mathbb R^N$ and its {\it support function\/}
$S(z):=\sup\{(x,z)\mid  x\in S\}$ for $z\in \mathbb R^N$.
Considering the members of $\mathbb R^N$ as singletons, we assume that
$\mathbb R^N$ lies in the set $\mathscr V_N$
of all compact convex subsets
of $\mathbb R^N$.

The Minkowski duality makes $\mathscr V_N$ into a~cone
in the space $C(S_{N-1})$  of continuous functions on the Euclidean unit sphere
$S_{N-1}$, the boundary of the unit ball $\mathfrak z_N$.
The
{\it linear span\/}
$[\mathscr V_N]$ of~$\mathscr V_N$ is dense in $C(S_{N-1})$, bears
a~natural structure of a~vector lattice
and is usually referred to as the {\it space of convex sets}.

The study of this space stems from the pioneering breakthrough of
Alexandrov in 1937 and the further insights of
Radstr\"{o}m, H\"{o}rmander, and Pinsker.

\section*{Linear Inequalities over Convex Surfaces}


 A measure $\mu$ {\it linearly majorizes\/} or {\it dominates\/}
a~measure $\nu$  on $S_{N-1}$ provided that to each decomposition of
$S_{N-1}$ into finitely many disjoint Borel sets $U_1,\dots,U_m$
there are measures $\mu_1,\dots,\mu_m$ with sum $\mu$
such that every difference $\mu_k - \nu|_{U_k}$
annihilates
all restrictions to $S_{N-1}$ of linear functionals over
$\mathbb R^N$. In symbols, we write $\mu\,{\gg}{}_{\mathbb R^N} \nu$.

Reshetnyak    proved  in 1954 (cp.~\cite{Reshetnyak}) that
$$
\int_{S_N-1} p d\mu \ge  \int_{S_N-1} p d\nu
$$
for each  sublinear    functional  $p$
on  $\mathbb R^N$   if   $\mu\,{\gg}{}_{\mathbb R^N} \nu$.
This gave an important trick for generating positive linear functionals
over various classes of convex  surfaces and functions.

\section*{Choquet's Order}


A~measure $\mu$ {\it affinely majorizes\/} or {\it dominates\/}
a measure $\nu$, both given     on a compact convex subset $Q$ of a locally convex space $X$,
provided that     to each decomposition of
$\nu$ into finitely many summands
$\nu_1,\dots,\nu_m$  there are measures $\mu_1,\dots,\mu_m$
whose sum is $\mu$ and for which every difference
$\mu_k - \nu_k$ annihilates all restrictions
to  $Q$  of affine   functionals over $X$.
In symbols, $\mu\,{\gg}{}_{\Aff(Q)} \nu$.

Cartier, Fell, and Meyer  proved in 1964  (cp.~\cite{Cartier_et_al}) that
$$
\int_{Q} f d\mu \ge  \int_{Q} f d\nu
$$
for each continuous convex function  $f$
on  $Q$   if and only if   $\mu\,{\gg}{}_{\Aff(Q)} \nu$.
An analogous necessity part for linear majorization was published
in 1969 (cp.~\cite{Kut69}--\cite{Dinges}).

\section*{Decomposition Theorem}
Majorization is a vast subject (cp.~\cite{Marshall_Olkin}).
The general form for many cones is as follows (cp.~\cite{Kut75}):


{\sl
Assume that $H_1,\dots,H_N$ are cones in a vector lattice~$X$.
Assume further that $f$ and $g$ are positive linear functionals on~$X$.
The inequality
$$
f(h_1\vee\dots\vee h_N)\ge g(h_1\vee\dots\vee h_N)
$$
holds for all
$h_k\in H_k$ $(k:=1,\dots,N)$
if and only if to each decomposition
of~$g$ into a~sum of~$N$ positive terms
$g=g_1+\dots+g_N$
there is a decomposition of~$f$ into a~sum of~$N$
positive terms $f=f_1+\dots+f_N$
such that
$$
f_k(h_k)\ge g_k(h_k)\quad
(h_k\in H_k;\ k:=1,\dots,N).
$$
}

\section*{Alexandrov Measures}

   Alexandrov proved the unique existence of
a translate of a convex body given its surface area function, thus completing the solution of
the Minkowski problem.
Each surface area function is an {\it Alexandrov measure}.
So we call a positive measure on the unit sphere which is supported by
no great hypersphere and which annihilates
singletons.

Each Alexandrov measure is a translation-inva\-riant
additive functional over the cone
$\mathscr V_N$.
The cone of positive translation-invariant measures in the
dual $C'(S_{N-1})$ of
 $C(S_{N-1})$ is denoted by~$\mathscr A_N$.

\section*{Blaschke's Sum}

 Given $\mathfrak x, \mathfrak y\in \mathscr V_N$,  the record
$\mathfrak x\,{=}{}_{\mathbb R^N}\mathfrak y$ means that $\mathfrak x$
and $\mathfrak y$  are  equal up to translation or, in other words,
are translates of one another.
So, ${=}{}_{\mathbb R^N}$ is the associate equivalence of
the preorder ${\ge}{}_{\mathbb R^N}$ on $\mathscr V_N$ of
the possibility of inserting one figure into the other
by translation.

The sum of the surface area measures of
$\mathfrak x$ and $\mathfrak y$ generates the unique class
$\mathfrak x\# \mathfrak y$ of translates which is referred to as the
{\it Blaschke sum\/} of $\mathfrak x$ and~$\mathfrak y$.

There is no need in  discriminating between a  convex figure,
the coset of its translates in  $\mathscr V_N/\mathbb R^N$,
and the corresponding measure in $\mathscr A_N$.

\section*{Comparison Between the Structures}

\begin{tabular}{|r|r|r|}
\hline
{\scshape Objects}&{\scshape Minkowski's Structure} &{\scshape
Blaschke's Structure}\\
\hline
cone of sets  &${\mathscr V}_N/\mathbb R^N$\hfil   &${\mathscr A}_N$\\
dual cone     &${\mathscr V}^*_N$\hfil             &${\mathscr A}^*_N$\hfil\\
positive cone &$\mathscr A^*_N $\hfil              &$\mathscr A_N$\hfil\\
linear functional &$V_1 (\mathfrak z_N,\,\cdot\,)$, breadth& $V_1(\,\cdot\,,\mathfrak z_N)$, area\hfil\\
concave functional &$V^{1/N}(\,\cdot\,)$\hfil           &$V^{(N-1)/N}(\,\cdot\,)$\hfil\\
convex program     & isoperimetric problem\hfil  & Urysohn's problem\hfil \\
operator constraint  & inclusion-like\hfil &  curvature-like\hfil \\
Lagrange's multiplier  &  surface\hfil & function\hfil \\
gradient\hfil        &$V_1(\bar{\mathfrak x},\,\cdot\,) $\hfil& $V_1(\,\cdot\,,\bar{\mathfrak x})$\hfil\\
\hline
\end{tabular}

\section*{The Natural Duality}

Let $C(S_{N-1})/\mathbb R^N$ stand for the factor space of
$C(S_{N-1})$ by the subspace of all restrictions of linear
functionals on $\mathbb R^N$ to $S_{N-1}$.
Let $[\mathscr A_N]$ be the space $\mathscr A_N-\mathscr A_N$
of translation-invariant measures, in fact, the linear span
of the set of Alexandrov measures.

$C(S_{N-1})/\mathbb R^N$ and $[\mathscr A_N]$ are made dual
by the canonical bilinear form
$$
\gathered
\langle f,\mu\rangle=\frac{1}{N}\int\nolimits_{S_{N-1}}fd\mu\\
(f\in C(S_{N-1})/\mathbb R^N,\ \mu \in[\mathscr A_N]).
\endgathered
$$

For $\mathfrak x\in\mathscr V_N/\mathbb R^N$ and $\mathfrak y\in\mathscr A_N$,
the quantity
$\langle {\mathfrak x},{\mathfrak y}\rangle$ coincides with the
{\it mixed volume\/}
$V_1 (\mathfrak y,\mathfrak x)$.

\section*{Solution of Minkowski's Problem}

Alexandrov observed that the gradient of $V(\cdot)$ at $\mathfrak x$ is proportional
to $\mu(\mathfrak x)$ and so minimizing $\langle \cdot,\mu\rangle$ over $\{V=1\}$
will yield the equality $\mu=\mu(\mathfrak x)$ by the  Lagrange multiplier
rule. But this idea fails since the interior of ${\mathscr V}_{N}$ is empty.
The fact that DC-functions are dense in $C(S_{N-1})$ is not helpful at all.

Alexandrov extended the volume to the positive cone of $C(S_{N-1})$ by the formula
$V(f):=\langle f,\mu(\co(f))\rangle$ with $co(f)$ the envelope of support functions
below $f$. The ingenious  trick settled all for the Minkowski problem.
This was done in 1938 but still is one of the summits of convexity.

In fact, Alexandrov suggested a functional analytical approach to extremal problems
for convex surfaces. To follow it directly in the general setting  is impossible
without the above description of the polar cones. The obvious limitations of the Lagrange multiplier rule are immaterial in the case of convex programs. It should be emphasized that the classical isoperimetric problem is not a Minkowski convex program in dimensions greater than~2. The convex counterpart is the Urysohn problem of maximizing volume given integral breadth \cite{Urysohn}.
The constraints of inclusion type are convex in the Minkowski structure,  which
opens way to complete solution  of new classes of Urysohn-type problems (cp.~\cite{Kut07}).

\section*{The External Urysohn Problem}
Among the convex figures, circumscribing $\mathfrak x_0 $ and having
integral breadth fixed, find a convex body of greatest volume.

{\sl
A feasible convex body $\bar {\mathfrak x}$ is a solution
to~the external Urysohn problem
if and only if there are a positive  measure~$\mu $
and a positive real $\bar \alpha \in \mathbb R_+$ satisfying

$(1)$ $\bar \alpha \mu
(\mathfrak z_N)\,{\gg}{}_{\mathbb R^N}\mu (\bar {\mathfrak x})+\mu $;

$(2)$~$V(\bar {\mathfrak x})+\frac{1}{N}\int\nolimits_{S_{N-1}}
\bar {\mathfrak x}d\mu =\bar \alpha V_1 (\mathfrak z_N,\bar {\mathfrak x})$;

$(3)$~$\bar{\mathfrak x}(z)={\mathfrak x}_0 (z)$
for all $z$ in the support of~$\mu $.
}

\section*{Solutions}
  If ${\mathfrak x}_0 ={\mathfrak z}_{N-1}$ then $\bar{\mathfrak x}$
is a {\it spherical lens} and $\mu$ is the restriction
of the surface area function
of the ball of radius
$\bar \alpha ^{1/(N-1)}$
to the complement of the support of the lens to~$S_{N-1}$.

If ${\mathfrak x}_0$ is an equilateral triangle then the solution
$\bar {\mathfrak x}$ looks as follows:

\centerline{\epsfxsize4cm\epsfbox{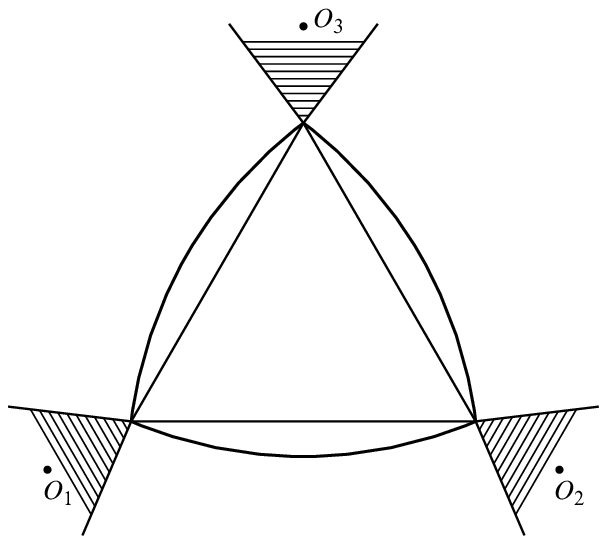}}

$\bar {\mathfrak x}$ is the union of~${\mathfrak x}_0$
and  three congruent slices of a circle of radius~$\bar \alpha$ and
centers $O_1$--$O_3$, while
$\mu$ is the restriction of $\mu(\mathfrak z_2)$
to the subset of $S_1$ comprising the endpoints
of the unit vectors of the shaded zone.

\section*{Symmetric Solutions}
This is  the general solution of the internal Urysohn problem inside a triangle
in the class of centrally symmetric convex figures:

\centerline{\epsfxsize6cm\epsfbox{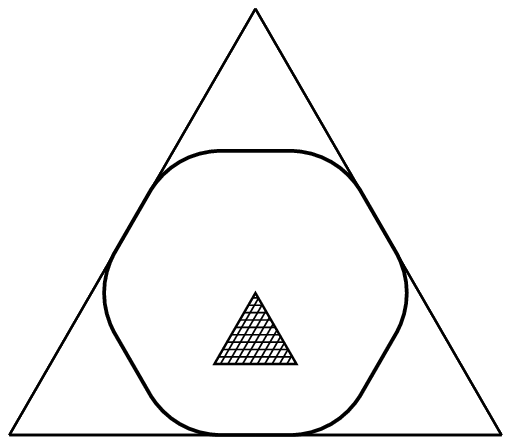}}

\section*{Current Hyperplanes}

  Find two convex figures $\bar{\mathfrak x}$ and $\bar{\mathfrak y}$
lying in a given convex body
$\mathfrak x_o$,
 separated by a~hyperplane with the unit outer normal~$z_0$,
and having the greatest total volume
of $\bar{\mathfrak x}$ and~$\bar{\mathfrak y}$
given the sum of their integral breadths.

{\sl
A feasible pair of convex bodies $\bar{\mathfrak x}$ and $\bar{\mathfrak y}$
solves the internal Urysohn problem with a current hyperplane
if and only if
there are convex figures  $\mathfrak x$ and $\mathfrak y$
and positive reals
$\bar\alpha $ and $\bar\beta$  satisfying

{\rm(1)}  $\bar{\mathfrak x}=\mathfrak x \# \bar\alpha\mathfrak z_N$;

{\rm(2)} $\bar{\mathfrak y}=\mathfrak y \# \bar\alpha\mathfrak z_N$;

{\rm(3)} $\mu(\mathfrak x) \ge \bar\beta\varepsilon_{z_0} $, $\mu(\mathfrak y) \ge \bar\beta\varepsilon_{-z_0} $;

{\rm(4)} $\bar {\mathfrak x}(z)=\mathfrak x_0 (z)$ for all $z\in \spt(\mathfrak x)\setminus \{z_0\} $;

{\rm(5)} $\bar {\mathfrak y}(z)=\mathfrak x_0 (z)$ for all $z\in \spt(\mathfrak x)\setminus \{-z_0\} $,
\noindent
with $\spt(\mathfrak x)$ standing for the {\it support\/} of $\mathfrak x$,
i.e. the support of the surface area measure $\mu(\mathfrak x)$
of~$\mathfrak x$.
}

\section*{Is Dido's Problem Solved?}

  From a utilitarian standpoint, the answer is
definitely in the affirmative. There is no evidence that Dido
experienced any difficulties, showed indecisiveness, and procrastinated the choice of the tract of land.
Practically speaking, the situation in which Dido  made her decision
was not as primitive as it seems at the first glance.

Assume that Dido
had known the isoperimetric property of the circle
and had been aware of the symmetrization processes that were elaborated
in the nineteenth century. Would this knowledge be sufficient for Dido
to choose the tract of land? Definitely, it would not.
The real coastline may  be rather ragged and craggy.
The photo snaps of coastlines are exhibited as the most visual
examples of fractality. From a theoretical standpoint, the free boundary
in Dido's planar problem  may be nonrectifiable, and so the concept of area
as the quantity to be optimized is itself rather ambiguous.
Practically speaking, the situation in which Dido  made her decision
was not as primitive as it seems at the first glance.
Choosing the tract of land, Dido had no right to trespass
the territory under the control of the local sovereign.
She had to choose the tract so as to encompass the camps of her
subjects and satisfy some fortification requirements.
Clearly, this generality is unavailable in the
mathematical models known as the classical isoperimetric problem.

Nowadays there is much research aiming at the problems with conflicting
goals (cp., for instance, \cite{MCDM}). One of the simplest and most
popular approach is based on the concept of Pareto-optimum.

\section*{ Pareto Optimality}
  Consider a~bunch of  economic agents
each of which intends to maximize his own income.
The {\it Pareto efficiency principle\/}  asserts
that  as an effective agreement of the conflicting goals it is reasonable
to take any state in which nobody can increase his income in any way other
than diminishing the income of at least one of the other fellow members.
 Formally speaking, this implies the search of the maximal elements
of the set comprising the tuples of incomes of the agents
at every state; i.e., some vectors of a finite-dimensional
arithmetic space endowed with the coordinatewise order. Clearly,
the concept of Pareto optimality was already abstracted to arbitrary
ordered vector spaces.

By way of example, consider a few multiple criteria problems of isoperimetric type.
For more detail, see \cite{Kut09}.

\section*{ Vector Isoperimetric Problem}

 Given are some convex bodies
$\mathfrak y_1,\dots,\mathfrak y_M$.
Find a convex body $\mathfrak x$ encompassing a given volume
and minimizing each of the mixed volumes $V_1(\mathfrak x,\mathfrak y_1),\dots,V_1(\mathfrak x,\mathfrak y_M)$.
In symbols,
$$
\mathfrak x\in\mathscr A_N;\
\widehat p(\mathfrak x)\ge \widehat p(\bar{\mathfrak x});\
(\langle\mathfrak y_1,\mathfrak x\rangle,\dots,\langle\mathfrak y_M,\mathfrak x\rangle)\rightarrow\inf\!.
$$
Clearly, this is a~Slater regular convex program in the Blaschke structure.

 {\sl
Each Pareto-optimal solution $\bar{\mathfrak x}$ of the vector isoperimetric problem
has the form}
$$
\bar{\mathfrak x}=\alpha_1{\mathfrak y}_1+\dots+\alpha_m{\mathfrak y}_m,
$$
where $\alpha_1,\dots,\alpha_m$ are positive reals.

\section*{The Leidenfrost Problem}

Given the volume of a three-dimensional
convex figure, minimize its  surface area and vertical  breadth.

By symmetry everything reduces to an analogous plane two-objective problem,
whose every Pareto-optimal solution is by~2 a~{\it stadium\/},
a weighted Minkowski sum of a disk and
a horizontal straight line segment.

{\sl A plane spheroid, a Pareto-optimal solution of the Leidenfrost problem,
is the result of rotation of a stadium around the vertical axis through
the center of the stadium}.

\section*{ Internal Urysohn Problem with Flattening}

Given are some~convex body
$\mathfrak x_0\in\mathscr V_N$ and some flattening direction~ $\bar z\in S_{N-1}$.
Considering $\mathfrak x\subset\mathfrak x_0$ of
fixed integral breadth, maximize the volume of~$\mathfrak x$ and  minimize the
breadth of $\mathfrak x$ in the flattening direction:
$\mathfrak x\in\mathscr V_N;\
\mathfrak x\subset{\mathfrak x}_0;\
\langle \mathfrak x,{\mathfrak z}_N\rangle \ge \langle\bar{\mathfrak x},{\mathfrak z}_N\rangle;\
(-p(\mathfrak x), b_{\bar z}(\mathfrak x)) \to\inf\!.
$

{\sl For a feasible convex body $\bar{\mathfrak x}$ to be Pareto-optimal in
the internal Urysohn problem with the flattening
direction~$\bar z$ it is necessary and sufficient that there be
positive reals $\alpha, \beta$ and a~convex figure $\mathfrak x$ satisfying}
$$
\gathered
\mu(\bar{\mathfrak x})=\mu(\mathfrak x)+ \alpha\mu({\mathfrak z}_N)+\beta(\varepsilon_{\bar z}+\varepsilon_{-\bar z});\\
\bar{\mathfrak x}(z)={\mathfrak x}_0(z)\quad (z\in\spt(\mu(\mathfrak x)).
\endgathered
$$

\section*{Rotational Symmetry}

Assume that a plane convex figure ${\mathfrak x}_0\in\mathscr V_2$ has the symmetry axis $A_{\bar z}$
with generator~$\bar z$.  Assume further that ${\mathfrak x}_{00}$ is the result of rotating
$\mathfrak x_0$  around the symmetry axis $A_{\bar z}$ in~$\mathbb R^3$.
$$
\gathered
\mathfrak x\in\mathscr V_3;\\
\mathfrak x  \text{\ is\ a\ convex\ body\ of\ rotation\ around}\ A_{\bar z};\\
\mathfrak x\supset{\mathfrak x}_{00};\
\langle {\mathfrak z}_N, \mathfrak x\rangle \ge \langle{\mathfrak z}_N,\bar{\mathfrak x}\rangle;\\
(-p(\mathfrak x), b_{\bar z}(\mathfrak x)) \to\inf\!.
\endgathered
$$

{\sl Each Pareto-optimal solution  is the result
of rotating around the symmetry axis a Pareto-optimal solution of the plane internal
Urysohn problem with flattening in the direction of the axis}.

\section*{Soap Bubbles}

Little is known about the analogous problems in arbitrary dimensions.
An especial place
is occupied by the result of Porogelov (cp.~ who
demonstrated that the ``soap bubble'' in a tetrahedron
has the form of the result of the rolling of a ball over a~solution
of the internal Urysohn problem, i.~e. the weighted Blaschke sum of
a tetrahedron and a ball.

\section*{The External Urysohn Problem with Flattening}

  Given are some convex body
$\mathfrak x_0\in\mathscr V_N$ and~flattening direction~$\bar z\in S_{N-1}$.
Considering $x\supset {\mathfrak x}_0$ of fixed integral breadth,
maximize volume and minimizing  breadth in
the flattening direction:
$
\mathfrak x\in\mathscr V_N;\
\mathfrak x\supset{\mathfrak x}_0;\
\langle \mathfrak x,{\mathfrak z}_N\rangle \ge \langle\bar{\mathfrak x},{\mathfrak z}_N\rangle;\
(-p(\mathfrak x), b_{\bar z}(\mathfrak x)) \to\inf\!.
$

{\sl For a feasible convex body $\bar{\mathfrak x}$ to be a Pareto-optimal solution of
the external Urysohn problem with flattening it is necessary and
sufficient that there be
positive reals $\alpha, \beta$, and a convex figure  $\mathfrak x$ satisfying}
$$
\gathered
\mu(\bar{\mathfrak x})+\mu(\mathfrak x)\gg{}_{{\mathbb R}^N} \alpha\mu({\mathfrak z}_N)+\beta(\varepsilon_{\bar z}+\varepsilon_{-\bar z});\\
V(\bar{\mathfrak x})+V_1(\mathfrak x,\bar{\mathfrak x})=\alpha V_1({\mathfrak z}_N,\bar{\mathfrak x})+ 2N\beta b_{\bar z}(\bar{\mathfrak x});\\
\bar{\mathfrak x}(z)={\mathfrak x}_0(z)\quad (z\in\spt(\mu(\mathfrak x)).
\endgathered
$$

\section*{ Optimal Convex Hulls}

Given  ${\mathfrak y}_1,\dots,{\mathfrak y}_m$ in~ $\mathbb R^N$,
place  ${\mathfrak x}_k$ within~ ${\mathfrak y}_k$, for $k:=1,\dots,m$,
maximizing the volume of
each of the  $\mathfrak x_1,\dots,\mathfrak x_m$ and minimize
the integral breadth of their convex hull:
$$
\mathfrak x_k\subset\mathfrak y_k ;\
(-p({\mathfrak x}_1),\dots,-p({\mathfrak x}_m), \langle \co\{{\mathfrak x}_1,\dots,{\mathfrak x}_m\},{\mathfrak z}_N\rangle)\to\inf.
$$

{\sl For  some feasible
${\bar{\mathfrak x}}_1,\dots,{\bar{\mathfrak x}}_m$  to have a
 Pareto-optimal convex hull it is necessary and sufficient
 that there be  $\alpha_1,\dots,\alpha_m\mathbb R_+$ not
 vanishing simultaneously and
 positive Borel  measures $\mu_1,\dots,\mu_m$ and $\nu_1,\dots, \nu_m$ on~$S_{N-1}$ such that}
 $$
\gathered
\nu_1+\dots+\nu_m=\mu({\mathfrak z}_N);\\
 \bar{\mathfrak x}_k(z)={\mathfrak y}_k(z)\quad (z\in\spt(\mu_k));\quad\\
\alpha_k \mu(\bar{\mathfrak x}_k)=\mu_k+\nu_k\ (k:=1,\dots,m).
\endgathered
$$


\end{document}